\documentclass[11pt]{article}
\usepackage{amstex}
\usepackage{amssymb}
\font\tenmsam=msam10
\def\square{\hbox{\tenmsam\char'003}}

\def\hfl#1#2{\smash{\mathop{\hbox to 12mm{\rightarrowfill}}
\limits^{\scriptstyle#1}_{\scriptstyle#2}}}

\def\mocat{$\mathcal C = (\mathcal C, \otimes, I)$ }
\def\socat{$\mathcal C = (\mathcal C, \otimes, I, \varphi)$ }
\def\iso{\widetilde {\longrightarrow}}

\newtheorem{theo}{Theorem}[section]
\newtheorem{prop}[theo]{Proposition}
\newtheorem{coro}[theo]{Corollary}
\newtheorem{lemm}[theo]{Lemma}
\newtheorem{defi}[theo]{Definition}
\newtheorem{rem}[theo]{Remark}
\newtheorem{ex}[theo]{Example}

\title{\bf Cosovereign Hopf algebras}
\date{Département des Sciences Mathématiques, Université Montpellier
II, \\ 
Place Eugène Bataillon, 34095 Montpellier Cedex 5, France} 
\author{Julien Bichon}

\makeatletter
\renewcommand{\@makefnmark}{}
\makeatother

\begin{document}
\maketitle

\hrule

\begin{abstract}
A sovereign monoidal category is an autonomous monoidal category
endowed with the choice of an autonomous structure and an
isomorphism of monoidal functors between the associated left and 
right duality functors. In this paper we define and study
the algebraic counterpart of sovereign monoidal categories:
cosovereign Hopf algebras.
In this framework we find a categorical characterization
of involutory Hopf algebras. We describe the universal
cosovereign Hopf algebras and we also study finite-dimensional
cosovereign Hopf algebras via the dimension theory provided
by the sovereign structure.
\end{abstract}

\smallskip


\smallskip

\hrule

\footnote{E-mail address: bichon\char64math.univ-montp2.fr}

\section{Introduction}

Monoidal category theory played a central role in the discovery of new
invariants of knots and links and in the development of the theory of 
quantum groups.

Let us recall that a tortile tensor category (or ribbon category)
is a braided monoidal category (\cite{[JS2]}), which is autonomous
(ie every object has a left dual, and hence also a right dual) and admits
a twist (\cite{[JS2], [Sh]}) compatible with duality.

The connection with knot theory is certainly best resumed is the following
coherence theorem by Shum (\cite{[Sh]}): `` {\sl the category of framed  
tangles (or tangles on ribbons) is the free tortile (or ribbon)
category generated by an object}''.
This means that to any object in a tortile (or ribbon) category, 
one can associate an isotopy invariant of framed tangles
(\cite{[Sh], [FY1], [Tu], [RT]}). We refer the reader to
the book \cite{[K]} for these topics.
\medskip

A new structure for monoidal categories appeared in papers
by Freyd and Yetter
(\cite{[FY2],[Ye]}). A sovereign structure on an autonomous
monoidal category (ie with left and right duals)
consists of the choice of a left and a right autonomous
structure and a monoidal isomorphism between the associated
left and right duality functors. A theorem of Deligne 
(proposition 2.11 in \cite{[Ye]}) brought interest 
in sovereign structures: there is a twist on 
an autonomous braided monoidal category if
and only if there is a sovereign structure on it.
 
Maltsiniotis (\cite{[Ma]}) studied sovereign monoidal categories in
their own right. He proposed a new equivalent definition
(which avoids the choice of an autonomous structure)
and he also showed that an axiom was redundant in \cite{[FY2],[Ye]}.
A sovereign structure is in fact the exact structure 
needed to define a trace theory.

\medskip

Quantum groups and monoidal categories are linked by tannakian duality
(\cite{[SR],[D],[U3]}): the reconstruction of a Hopf algebra
from its finite-dimensional comodules.
The use of tannakian duality is clearly illustrated in \cite{[JS1]}:
to an additional categorical structure on the category of finite-dimensional
comodules, one associates an additional 
algebraic structure on the Hopf algebra.
For example if $A$ is a Hopf algebra, then there is a braiding
(\cite{[JS2]}) on ${\rm Co}_f(A)$ (finite-dimensional $A$-comodules)
if and only if there is a cobraiding on $A$, ie a linear form on
$A \otimes A$ satisfying certain conditions (see the appendix or 
\cite{[K],[JS1]}). In the same way if $A$ is a cobraided 
Hopf algebra, the category ${\rm Co}_f(A)$ is balanced
(there is a twist on it) if and only if there is a linear form
$\tau$ (called a cotwist) on $A$ satisfying some conditions.
Another example of the use of Tannaka duality
is given by Street in \cite{[St]}: the construction
of the quantum double for any bialgebra (even infinite-dimensional).
\medskip

In this paper we find the algebraic structure on Hopf algebras
corresponding to sove-reign structures. A sovereign character
on a Hopf algebra $A$ with bijective antipode is a {\sl character}
$\Phi$ on $A$ such that $S^{-1} = \Phi * S * \Phi^{-1}$
($S$ is the antipode of $A$ and $*$ is the convolution product).
A cosovereign Hopf algebra is a pair $(A,\Phi)$ where $A$
is a Hopf algebra with bijective antipode and $\Phi$
is a sovereign character on $A$.

Let $A$ be a Hopf algebra with a bijective antipode.
We show that the category ${\rm Co}_f(A)$ admits a sovereign structure
if and only if there is a sovereign character on $A$
(in fact the main result, theorem 3.12, is more precise).
We also obtain a categorical characterization of involutory
Hopf algebras (Hopf algebras whose square of the antipode is equal to the 
identity): a Hopf algebra is involutory if and only if
the category ${\rm Co}_f(A)$ admits a sovereign structure for which 
the forgetful functor is sovereign.

\medskip

Back not too far from knot theory, the theorem of Deligne mentioned
above states in particular the bijective 
correspondence of cotwists and sovereign
characters for a cobraided Hopf algebra (we give a purely Hopf algebraic
proof of this result in an appendix). But it is easier to check
the existence of a sovereign character (since it is a character).
Therefore a technical simplification is brought by sovereign structures
in this context.  
     
\medskip

The end of the paper is devoted to the study of some examples.
We describe the universal (or free)
cosovereign Hopf algebras: every finite-type
cosovereign Hopf algebra is a homomorphic quotient of one of them.
These algebras are parameterized by an invertible matrix.
When the base field is the field of complex numbers, they already
appeared in a different context: they are the algebras of representative
functions on the universal compact quantum groups defined by Van Daele and
Wang (\cite{[VW]}). We also examine a class of examples closely related
to the quantum groups $SU(n)$ of Woronowicz (\cite{[W2]}) 
and show that these Hopf algebras are cosovereign.

We also study finite-dimensional cosovereign Hopf algebras via the dimension
theory provided  by the sovereign structure (theorem 5.1).
 We show that if $(A,\Phi)$
is a  finite-dimensional cosovereign Hopf algebra 
over a field of characteristic zero whose internal
dimension (ie the dimension computed in the sovereign monoidal
category ${\rm Co}_f(A)$) is non-zero, then $A$ is involutory
(and therefore is semisimple and cosemisimple by \cite{[La], [LR]}).

\medskip

The paper is organized as follows. In section 2 we review the
basic definitions of monoidal category theory.
In section 3 we introduce cosovereign (and sovereign)
Hopf algebras and study their relations with sovereign
structures for monoidal categories.
Sections 4 and 5 are devoted to examples.
In an appendix we examine the relations between sovereign characters
and cotwists for cobraided Hopf algebras.

\subsection*{Notations}

Throughout this paper $k$ will denote a commutative field.
The category of finite-dimen-sional vector spaces will be denoted by
${\rm Vect}_f(k)$.

We assume the reader to be familiar with the theory of Hopf algebras
(\cite{[A], [Sw], [K]}) and in particular we freely use
convolution products.

Let  $A = (A,m,u,\Delta,\varepsilon,S)$ be a Hopf $k$-algebra.
 The multiplication will be denoted by  $m$,   $u : k
\rightarrow A$ is the unit of $A$, while $\Delta$, $\varepsilon$ and $S$
are respectively the  comultiplication, the counit and the antipode of $A$.
 
If $A$ is a Hopf algebra, the category of finite-dimensional right 
$A$-comodules will be denoted by  ${\rm Co}_f(A)$
while the category of finite-dimensional left $A$-modules
will be denoted by  ${\rm Mod}_f(A)$. We only consider
right comodules and left modules.

\section{Monoidal categories}

The aim of this section is to recall the basic definitions of monoidal
category theory and to fix some notations. We refer the reader to
\cite{[ML]} or \cite{[JS2]} for the general definitions of monoidal categories and monoidal functors. The material presented here is now classical
and hence we will be a little concise.
We only use strict monoidal categories:
by Mac-Lane's coherence theorem (see \cite{[JS2]}, 1.4 for a simple proof), every monoidal category is monoidally equivalent to a strict one.
However the reconstruction theorem (2.12) we use deals with non strict
monoidal categories. We will see that the general statement follows from
the particular case of strict monoidal categories.

\bigskip

\begin{defi}
A monoidal category $\mathcal C = (\mathcal C, \otimes, I)$ consists of a 
category $\mathcal C$, a functor $\otimes : \mathcal C \times \mathcal C
\longrightarrow \mathcal C$ (called the tensor product of $\mathcal C$)
and an object $I$ of $\mathcal C$ (called the monoidal unit) such that 
for all objects $X,Y,Z$ and all arrows $f,g,h$ of $\mathcal C$ we have:

\noindent
1)
associativity : $(X \otimes Y)\otimes Z = X \otimes (Y \otimes Z)$ and
$(f \otimes g) \otimes h = f \otimes (g \otimes h)$.

\noindent
2) unit : $ I \otimes X = X = X \otimes I$ and 
$1_I \otimes f = f = f \otimes 1_I$.
\end{defi}

Let \mocat  be a monoidal category. Then $(\mathcal C^o, \otimes, I)$,
$(\mathcal C, \otimes^o, I)$ and $(\mathcal C^o, \otimes^o, I)$
are also monoidal categories, where $\mathcal C^o$ is the opposite
category of $\mathcal C$ and $\otimes^o$ is the opposite tensor product of
$\mathcal C$ ($X \otimes^o Y=Y \otimes X$). 

\begin{defi}
Let $\mathcal C$ and $\mathcal D$ be monoidal categories.
A monoidal functor $F=(F, \widetilde F)$ consists of a functor
$F : \mathcal C \longrightarrow \mathcal D$, a natural family
of isomorphisms $\widetilde F_{X,Y} : F(X) \otimes F(Y) \iso F(X \otimes Y)$ 
and an isomorphism $\widetilde F_I :I \iso F(I)$ such that:

$\widetilde F_{X \otimes Y,Z} \circ (\widetilde F_{X,Y} \otimes 1_{F(Z)})
= \widetilde F_{X, Y \otimes Z} \circ (1_{F(X)} \otimes \widetilde F_{Y,Z}$),

$\widetilde F_{I,X} \circ (\widetilde F_I \otimes 1_{F(X)}) = 1_{F(X)} =
\widetilde F_{X,I} \circ (1_{F(X)} \otimes \widetilde F_I)$.

A monoidal equivalence is a monoidal functor whose underlying functor is an 
equivalence of categories.

Let $F, G : \mathcal C \longrightarrow \mathcal D$ be monoidal functors.
A morphism of monoidal functors $u : F \longrightarrow G$ is a natural
transformation $u : F \longrightarrow G$ such that:

$u_{X \otimes Y} \circ \widetilde F_{X,Y} = \widetilde G_{X,Y} \circ
(u_X \otimes u_Y)$ and $u_I \circ \widetilde F_I = \widetilde G_I$.

If $u$ is also an isomorphism we write $F \cong^{\otimes} G$.

\noindent
The set of morphisms between monoidal functors $F$ and $G$ will be
denoted by ${\rm Hom}^{\otimes}(F,G)$.
\end{defi}

Let $F : \mathcal C \longrightarrow \mathcal D$ be a monoidal equivalence.
By \cite{[SR]} there is a monoidal equivalence 
$G : \mathcal D \longrightarrow \mathcal C$ and isomorphisms of
monoidal functors $1_{\mathcal C} \cong^\otimes G \circ F$ and
$1_{\mathcal D} \cong^\otimes F \circ G$.

Let $F, G : \mathcal C \longrightarrow \mathcal D$ be monoidal functors
and let $u : F \longrightarrow G$ be a morphism of monoidal functors.
Let $K : \mathcal D \longrightarrow \mathcal E$ be a monoidal functor. Then
$K(u) : KF \longrightarrow KG$ is a morphism of monoidal functors.
If $u$ is an isomorphism, so is $K(u)$. Let $K' : \mathcal D \longrightarrow \mathcal E$ be another monoidal functor and let $v: K \longrightarrow K'$
be a morphism of monoidal functors. Then $v_F : KF \longrightarrow K'F$
is a morphism of monoidal functors. If $v$ is an isomorphism, so is $v_F$.

\begin{defi}
Let \mocat be a monoidal category and let $X \in {\rm ob}(\mathcal C)$.
A left dual for X is a triplet $( {^\vee X},\varepsilon_X,\eta_X)$ with
$^\vee X \in {\rm ob}(\mathcal C)$, 
$\varepsilon_X :  {^\vee X} \otimes X \longrightarrow I$ and 
$\eta_X : I \longrightarrow X \otimes  {^\vee X}$ are morphisms of $\mathcal C$
such that:
$$ 
(1_X \otimes \varepsilon_X) \circ (\eta_X \otimes 1_X) = 1_X \quad 
{\rm and} \quad (\varepsilon_X \otimes 1_{^\vee X})
 \circ (1_{^\vee X} \otimes \eta_X) = 1_{^\vee X}.
$$
A right dual for $X$ is a triplet $(X^\vee, e_X, d_X)$
with $X^\vee \in {\rm ob}(\mathcal C)$,
$e_X :X \otimes X^\vee \longrightarrow I$ and 
$d_X : I \longrightarrow X^\vee \otimes X$ are morphisms of $\mathcal C$
such that:
$$ 
(1_{X^\vee} \otimes e_X) \circ (d_X \otimes 1_{X^\vee}) = 1_{X^\vee} \quad 
{\rm and} \quad (e_X \otimes 1_X)
 \circ (1_X \otimes d_X) = 1_X.
$$
\end{defi}

Let $X$ be an object of a monoidal category and suppose that $X$ is endowed
with a left dual. Then the functor $X \otimes -$ admits a left adjoint
$^\vee X \otimes -$ (see \cite{[JS1]}, Section 9). Similarly,
if there is a right dual for $X$, the functor $X \otimes -$ admits
a right adjoint which is $X^\vee \otimes -$.
We inherit the unicity results of adjoint functors:

\begin{prop}
Let $X$ be an object in a monoidal category $\mathcal C$. Suppose that
$(^\vee X,\varepsilon_X,\eta_X)$ and $(^*X,\varepsilon_X',\eta_X')$
are left duals for $X$. Then there is a unique isomorphism
$l_X :  {^\vee X} \longrightarrow  {^*X}$ such that:
$$
\varepsilon_X' \circ (l_X \otimes 1_X) = \varepsilon_X \quad
({\rm and} \quad (1_X \otimes l_X) \circ \eta_X = \eta'_X).
$$
Suppose that $(X^\vee, e_X, d_X)$ and $(X^*, e_X', d_X')$ are
right duals for $X$. Then there is a unique isomorphism
$r_X : X^\vee \longrightarrow X^*$ such that:
$$
e_X' \circ (1_X \otimes r_X) =e _X \quad
({\rm and} \quad (r_X \otimes 1_X) \circ d_X = d'_X).
$$
\end{prop}

\begin{defi}
Let $\mathcal C$ be a monoidal category. Let $X,Y \in {\rm ob}(\mathcal C)$
and suppose that that $X$ and $Y$ admit a left dual (resp. a right dual).
Let $f \in {\rm Hom}_{\mathcal C}(X,Y)$. The left transpose of $f$
(resp. right transpose of $f$) is the morphism
$^t f :  {^\vee Y}
 \longrightarrow  ^\vee $X (resp. $f^t : Y^\vee \longrightarrow
X^\vee$ defined by:
$$^t f = (\varepsilon_Y \otimes 1_{^\vee X}) \circ
(1_{^\vee Y} \otimes f \otimes 1_{^\vee X}) \circ
(1_{^\vee Y} \otimes \eta_X) 
$$
$$({\rm resp.} \ \
f^t = (1_{X^\vee}\otimes e_Y) \circ
(1_{X \vee} \otimes f \otimes 1_{Y^\vee}) \circ
(d_Y \otimes 1_{Y^\vee})).
$$
\end{defi}

It is easily seen that $^tf$ (resp. $f^t$
is the only morphism satisfying
$\varepsilon_X \circ (^tf \otimes 1_X) 
= \varepsilon_Y \circ (1_{^\vee Y} \otimes f)$ or 
$(1_Y \otimes ^tf) \circ \eta_Y = (f \otimes 1_{^\vee X}) \circ \eta_X$
(resp.  
$e_X \circ (1_X \otimes f^t) 
= e_Y \circ (1_{^\vee Y} \otimes f)$ or 
$(f^t \otimes 1_Y) \circ d_Y = (1_{X^\vee} \otimes f) \circ d_X).$

\begin{defi}
A monoidal category is said to be left autonomous (resp. right autono-mous ;
autonomous) if every object has a left dual (resp. right dual ; resp left
and right duals). 
\end{defi}
 
\begin{defi}
Let $\mathcal C$ be a left (resp.right) autonomous monoidal category.
A left (resp. right) autonomous structure on $\mathcal C$ is the choice
 for every object of a left dual $({^\vee X},\varepsilon_X,\eta_X)$
(resp. a right dual  $(X^\vee, e_X, d_X)$) such that $^\vee I =I$
and $\varepsilon_I = \eta_I = 1_I$ (resp. $I^\vee = I$ and $e_I = d_I= 1_I$).
An autonomous structure on an autonomous monoidal category $\mathcal C$
consists of a left autonomous structure on $\mathcal C$ and a right
autonomous structure on $\mathcal C$.
\end{defi}

\begin{defi}
Let \mocat be an autonomous monoidal category.
Choose an autonomous structure on $\mathcal C$.
We get two monoidal functors: 

${\bf D}_l : (\mathcal C, \otimes, I) \longrightarrow 
(\mathcal C^o, \otimes^o, I)$ defined by ${\bf D}_l(X) =  {^\vee X}$ and
${\bf D}_l(f) = {^tf},$

${\bf D}_r : (\mathcal C, \otimes, I) \longrightarrow 
(\mathcal C^o, \otimes^o, I)$, defined by ${\bf D}_r(X) = X^\vee$ and
${\bf D}_r(f) = f^t,$

\noindent
called the left duality functor and the right duality functor respectively.
\end{defi}

Another choice of autonomous structure would lead to monoidal functors
${\bf D}_l'$ and ${\bf D}_r'$ with ${\bf D}_l' \cong^\otimes {\bf D}_l$
and ${\bf D}_r' \cong^\otimes {\bf D}_r$ (see proposition 2.4). Thus 
the choice of an autonomous structure is just a convenient way to define the
duality functors.

\medskip

Let us remark that the duality functors can also be seen as monoidal functors
$(\mathcal C^o, \otimes^o, I)\break \longrightarrow (\mathcal C, \otimes, I)$.
There are isomorphisms of monoidal functors (see proposition 2.4):
$$h : 1_{\mathcal C} \cong^\otimes {\bf D}_l \circ {\bf D}_r 
\quad (2.8.1) \quad
\kappa : 1_{\mathcal C} \cong^\otimes  {\bf D}_r \circ {\bf D}_l 
\quad (2.8.2)$$
Let $\mathcal C$ and $\mathcal D$ be autonomous monoidal categories
endowed with autonomous structures and let $F : \mathcal C \longrightarrow
\mathcal D$ be a monoidal functor. Then there are isomorphisms of monoidal
functors (see proposition 2.4 again):
$$l : {\bf D}_l \circ F \cong^\otimes F \circ {\bf D}_l 
\quad (2.8.3) \quad 
r : {\bf D}_r \circ F \cong^\otimes F \circ {\bf D}_r 
\quad (2.8.4)$$

\begin{ex}
{\rm We briefly describe the autonomous structure on
${\rm Vect}_f(k)$ and ${\rm Co}_f(A)$ where $A$ a Hopf algebra with
bijective antipode. Let $V$ be a finite-dimensional vector space.
Let $^\vee V = V^\vee = V^* = {\rm Hom}(V,k)$.
It is well known that this procedure, with classical evaluation and coevaluation maps, defines an autonomous structure on ${\rm Vect}_f(k)$.
Now let $V$ be a finite-dimensional $A$-comodule with coaction
$\alpha_V : V \longrightarrow V \otimes A $ such that
$\alpha_V(v_i) = \sum_j v_j \otimes a_{ji}$ for some basis $(v_i)$ of V.
Now let $^\vee V$ (resp. $V^\vee$) be the $A$-comodule whose underlying
vector space is $V^*$ and whose coaction 
$\alpha_{^\vee V} :  {^\vee V} \longrightarrow ^\vee V \otimes A$
(resp. $\alpha_{V^\vee } : V^\vee \longrightarrow V^\vee  \otimes A$) 
is defined by
$\alpha_{^\vee V}(v_i^*) = \sum_j v_j^* \otimes S(a_{ij})$
(resp. $\alpha_{V^\vee}(v_i^*) = \sum_j v_j^* \otimes S^{-1}(a_{ij})$)
where $(v_i^*)$ is the dual basis of $(v_i)$.
Then $^\vee V$ and $V^\vee$, with classical evaluation and coevaluation maps,
are left and right duals for $V$.}
\end{ex}

\medskip

We now proceed to describe the tannakian reconstruction theorems.
First we need the following definition as given in \cite{[JS1]}:

\begin{defi}
Let $\mathcal C$ be a small category and let
$F : {\mathcal C} \longrightarrow {\rm Vect}_f(k)$ be a functor.
Let $\mathcal N$ be the subspace of
$\bigoplus_{X \in {\rm ob} \mathcal C} {\rm End}(F(X))$
generated by the expressions $g \circ F(f) - F(f) \circ g$,
where $f \in {\rm Hom}_{\mathcal C}(X,Y)$ 
and $g \in {\rm Hom}(F(Y),F(X))$. We define
$$ {\rm End}^\vee(F) = \bigoplus_{X \in {\rm ob} \mathcal C} {\rm End}(F(X))
/ \mathcal N$$
Let $g \in {\rm End}(F(X))$. We denote it by $[X,g]$ as an element of  
${\rm End}^\vee(F)$.
\end{defi}

\begin{theo}
Let $\mathcal C$ be a small category and let
$F : {\mathcal C} \longrightarrow {\rm Vect}_f(k)$ be a functor.

\noindent
i) The vector space ${\rm End}^\vee(F)$ is a coalgebra with coproduct given
by $$\Delta([X,\phi\otimes x])= \sum_i [X, \phi\otimes
v_i] \otimes [X, v_i^* \otimes x],$$ where $X$ is an object of $\mathcal C$,
$\phi \in F(X)^*$ , $x\in F(X)$ and 
$(v_i)$ is a basis of $F(X)$. The counit is given by 
$\varepsilon([X,g]) = {\rm Tr}(g)$. There is a linear isomorphism
$${\rm End}(F) \longrightarrow {\rm End}^\vee(F)^*$$ (where ${\rm End}(F)$ is the algebra of endomorphisms of the functor $F$) which to
$u \in {\rm End(F)}$, associates the linear form $f_u$ defined by
$f_u([X,g])= {\rm Tr}(u_X \circ g)$.

\smallskip

\noindent
ii) Suppose now that \mocat is a monoidal category and that 
$F = (F,\widetilde F)$ is a monoidal functor. Then there is 
a bialgebra structure on ${\rm End}^\vee(F)$  whose product is given by
$$
[X,g].[Y,h] = [X \otimes Y, \widetilde F_{X,Y} \circ (g \otimes h)
\circ \widetilde F_{X,Y}^{-1}]
$$
and whose unit is $[I,1]$.

\smallskip

\noindent
iii) Suppose that $\mathcal C$ is a left autonomous monoidal category
and that $F$ is a monoidal functor.
Then the bialgebra ${\rm End}^\vee(F)$ is a Hopf algebra with antipode 
$S$ given by $$S([X,g]) = [^\vee X,l_X \circ  ^tg \circ l_X^{-1}]$$
($l_X$ is from 2.8.3). The linear isomorphism of i) induces a bijection
$${\rm Aut}^\otimes (F) \longrightarrow {\rm Hom}_{k-{\rm alg}}({\rm End}^\vee(F),k)$$
where ${\rm Aut}^\otimes (F)$ denotes the set of endomorphisms of 
the monoidal functor $F$
(which are automorphisms). If furthermore $\mathcal C$ is autonomous,
then the antipode of ${\rm End}^\vee(F)$ is bijective, with inverse
given by $$S^{-1}([X,g]) = [X^\vee, r_X \circ g^t \circ r_X^{-1}])$$
($r_X$ is from 2.8.4).    
\end{theo}

The proof of this theorem can be found in \cite{[JS1]} for example.
Statements i) and ii) are due to Saavedra \cite{[SR]}.
Statement iii) is due to Ulbrich \cite{[U3]}. This theorem is the framework
for the following reconstruction theorem, whose  proof can be found
in many places and at different levels of generality (\cite{[SR], [D],
[Br], [JS1]}).

\begin{theo}
i) Let $\mathcal C$ be an autonomous monoidal (small) category and let
$F : {\mathcal C} \longrightarrow {\rm Vect}_f(k)$ be a monoidal functor.
Then $F$ factorizes through a monoidal functor
$\overline F : {\mathcal C} \longrightarrow {\rm Co}_f({\rm End}^\vee(F))$
followed by the forgetful functor.

\smallskip

\noindent
ii) Let $\mathcal C$ be a (small) $k$-tensorial autonomous monoidal category
and let $F : {\mathcal C} \longrightarrow {\rm Vect}_f(k)$ be a 
fibre functor. Then the above functor $\overline F$ is a monoidal
equivalence of categories.

\smallskip

\noindent
iii) Let $\mathcal C = {\rm Co}_f(A)$ be be the 
category of right $A$-comodules 
where $A$ is a Hopf algebra with bijective antipode and let $F$
be the forgetful functor. Then the Hopf algebras $A$ and  
${\rm End}^\vee(F)$ are isomorphic.
\end{theo}

We have just included statement ii) for the sake of completeness.
We refer the reader to \cite{[Br]} for the precise definitions.

\begin{rem}
{\rm In theorem 2.11 and 2.12, the encountered monoidal
categories must not be assumed to be strict for useful applications.
This is clear from statement iii) in theorem 2.12. However, assume that theorem
2.11 has been proved for strict monoidal categories. Let $\mathcal C$
be a (non-strict) monoidal category and let 
 $F : {\mathcal C} \longrightarrow {\rm Vect}_f(k)$ be a monoidal functor.
Let $St(\mathcal C)$ be a strict monoidal category with a
monoidal equivalence $i : St(\mathcal C) \longrightarrow \mathcal C$.
Then the coalgebras ${\rm End}^\vee(F)$ and ${\rm End}^\vee(F \circ i)$
are easily seen to be isomorphic and therefore statements ii) and iii)
of theorem 2.11 follow in the general case (and the same formulas hold);
In the same way, statements i) and ii) of theorem 2.12 follow from the 
strict monoidal case.

Let $A$ be a Hopf algebra. Assume that there is an additional
categorical structure on ${\rm Co}_f(A)$ : we would like to translate it
into an additional algebraic structure on $A$. It is clear from
the precedent discussion that it is sufficient to make this
translation for ${\rm End}^\vee(F)$, where 
$F : {\mathcal C} \longrightarrow {\rm Vect}_f(k)$ is a monoidal
functor and $\mathcal C$ is a strict monoidal category.}
\end{rem}

\begin{rem}
{\rm When $k$ is a ring theorem 2.11 is still valid in some sense:
one has to replace Vect$_f(k)$ by Proj$_f(k)$, the category of 
finitely generated projective $k$-modules.
Let us note that if $\mathcal C$ is autonomous, a monoidal functor
$F : \mathcal C \longrightarrow {\rm Mod}(k)$
take values in  Proj$_f(k)$ (\cite{[D]}, 2.6).}
\end{rem}

\section{Sovereign monoidal categories and cosovereign Hopf algebras} 

In this section we give our main result: the correspondence
between categorical and algebraic sovereign structures.

\begin{defi}
Let \mocat be an autonomous monoidal category. A sovereign structure
$\varphi$ on $\mathcal C$ consists of an autonomous structure on 
$\mathcal C$ and an isomorphism of monoidal functors 
$\varphi : {\bf D}_r \cong^{\otimes} {\bf D}_l$ for the duality functors.
A sovereign monoidal category \socat is an autonomous monoidal category
endowed with a sovereign structure.
\end{defi}

\begin{rem}
{\rm Let \socat be a sovereign monoidal category. Let us choose another 
autonomous structure on $\mathcal C$ and let us denote by
${\bf D}'_l$ and ${\bf D}'_r$ the associated duality functors.
Then there is an isomorphism ${\bf D}_r' \cong^\otimes
{\bf D}_r \ ^{\varphi} \cong^\otimes {\bf D}_l \cong ^\otimes {\bf D}_l'$.
In this way we get another sovereign structure $\varphi'$ on
$\mathcal C$. This fact suggests that we possibly could give a definition
of sovereign structure which does not depend on the choice of an
autonomous structure. This is done by Maltsiniotis in definition
3.1.2 of \cite{[Ma]}. He shows in theorem 3.2.2
that his definition coincides with definition 3.1 given here. 

In the earlier definition by Freyd and Yetter (\cite{[FY2], [Ye]}) there 
was another axiom, which was shown to be redundant by Maltsiniotis (proposition
3.2.3 in \cite{[Ma]}).} 
\end{rem}

\begin{defi} 
Let $(\mathcal C, \otimes, I, \varphi)$ and 
 $(\mathcal D, \otimes, I, \psi)$ be sovereign monoidal categories.
 A\break monoidal functor
$F=(F,\widetilde F) :  (\mathcal C, \otimes, I)
\longrightarrow  (\mathcal D, \otimes, I)$ is said to be sovereign if 
we have $F(\varphi) \circ r = l \circ \psi_F$ where $l$ and
$r$ are the isomorphisms of (2.8.3) and (2.8.4).   
\end{defi}

\begin{ex}
{\rm The category of diagrams (\cite{[Ye]}) is the free sovereign monoidal
category generated by an object (\cite{[FY2]}).

An autonomous braided monoidal category admits a sovereign structure
if and only there is a twist on it (\cite{[Ye]}, see also \cite{[Ma]})}.
\end{ex}

\medskip

The following lemmas will be useful :

\begin{lemm}
Let \socat and $(\mathcal D,\otimes,I,\psi)$ be sovereign monoidal
categories and let
$F=(F,\widetilde F) :  (\mathcal C, \otimes, I)
\longrightarrow  (\mathcal D, \otimes, I)$ be a monoidal functor.
There is an element $u^{\varphi,\psi} \in {\rm Aut}^\otimes(F)$ defined by
$$
u^{\varphi,\psi} = \kappa_F^{-1} \circ {\bf D}_r(l^{-1}
\circ F(\varphi) \circ r \circ \psi_F^{-1}) \circ \kappa_F
$$
(where $\kappa$, $l$ and $r$ are the isomorphisms of (2.8.2), (2.8.3)
and (2.8.4)). If furthermore $F$ is sovereign, then
$u^{\varphi,\psi} = 1_F$. 
\end{lemm}

\noindent
{\bf Proof}. The element $u^{\varphi,\psi}$ of $ {\rm End}(F)$ just defined
is the composition of isomorphisms of monoidal functors.
Hence $u^{\varphi,\psi}$ is an automorphism of the monoidal functor $F$.
If $F$ is a sovereign functor, then
$l^{-1} \circ F(\varphi) \circ r \circ \psi_F^{-1}$ is the identity
morphism of ${\bf D}_l \circ F$, and thus $u^{\varphi,\psi}=1_F$. \square

\begin{lemm}
Let \socat and 
$\mathcal D =(\mathcal D,\otimes,I,\psi)$ be sovereign monoidal
categories and let
$F : (\mathcal C, \otimes, I,\varphi)
\longrightarrow  (\mathcal D, \otimes, I,\phi)$ be 
a sovereign monoidal functor. Suppose that we have another choice of 
autonomous structure on $\mathcal D$ and let $\psi'$ the associated sovereign 
structure of remark 3.2. Then
$F :  (\mathcal C, \otimes, I, \varphi)
\longrightarrow  (\mathcal D, \otimes, I, \psi')$ is still a sovereign functor.
\end{lemm}

\noindent
{\bf Proof}. Let ${\bf D}_r'$ and ${\bf D}_l'$ be the duality functors on $\mathcal D$. Let $\alpha : {\bf D}_r' \cong^\otimes {\bf D}_r$ and
$\beta : {\bf D}_l' \cong^\otimes {\bf D}_l$ be the isomorphisms
given by proposition 2.4: $\psi' = \beta^{-1} \circ \psi \circ \alpha$.
Let $r' : {\bf D}_r'\circ F \cong^\otimes F \circ {\bf D}_r$,
$l' : {\bf D}_l'\circ F \cong^\otimes F \circ {\bf D}_l$,
$r : {\bf D}_r \circ F \cong^\otimes F \circ {\bf D}_r$,
$l : {\bf D}_l\circ F \cong^\otimes F \circ {\bf D}_l$
as in (2.8.3) and (2.8.4). Then we have $r' = r \circ \alpha_F$ and
$l' = l \circ \beta_F$ (since they verify the equations in proposition 2.4).
Thus we have $F(\varphi) \circ r' = F(\varphi) \circ r \circ \alpha_F
= l \circ \psi_F \circ \alpha_F = l' \circ \psi_F'$. \square
\medskip

We now come into the heart of the subject:

\begin{defi}
Let $A$ be a Hopf algebra with bijective antipode. A sovereign character on
$A$ is a character $\Phi$ on $A$ such that
$S^{-1}= \Phi * S * \Phi^{-1}$. A cosovereign Hopf algebra is a
pair $(A,\Phi)$ where $A$ is a Hopf algebra with bijective antipode 
and $\Phi$ is a sovereign character on $A$. 
\end{defi}

\begin{rem}
{\rm We could have defined a cosovereign Hopf algebra as a
pair $(A, \Phi)$ consisting of a Hopf algebra $A$ and a character $\Phi$ on $A$
such that $S^2 = \Phi^{-1} * id * \Phi$. The bijectivity of the antipode follows immediately. We have chosen this definition since it is more natural
from the categorical viewpoint.}
\end{rem}

There is an immediate dual definition:

\begin{defi}
Let $A$ be a Hopf algebra with bijective antipode. A sovereign element
of $A$ is a group-like element $\Phi$ such that
$S^{-1}(a) = \Phi S(a) \Phi^{-1}$ for all $a \in A$.
A sovereign Hopf algebra is a pair $(A,\Phi)$ where $A$ is a Hopf algebra with 
bijective antipode and $\Phi$ is a sovereign element of $A$.
\end{defi}

It is clear that if $(A,\Phi)$ is a cosovereign (resp. sovereign)
Hopf algebra then $(A^0,\Phi)$ is a sovereign (resp. cosovereign)
Hopf algebra.

\begin{prop}
Let $(A,\Phi)$ be a cosovereign Hopf algebra. Then the sovereign
character $\Phi$ defines a sovereign structure on ${\rm Co}_f(A)$.
\end{prop} 

\noindent
{\bf Proof}. We use the autonomous structure on ${\rm Co}_f(A)$
described in example 2.9. Let $V$ be a finite-dimensional
$A$-comodule with coaction $\alpha_V$. We define a linear map
$\varphi_V  : V^\vee \longrightarrow {^\vee V}$ as follows:
$\varphi_V = (1_{^\vee V} \otimes \Phi^{-1}) \circ \alpha_{^\vee V}$
(we use the fact $^\vee V = V^\vee$ as vector space).
For the convenience of the reader, let us describe $\varphi_V$ 
explicitly. If $(v_i)$ is a basis of $V$ such that
$\alpha_V(v_i) = \sum_j v_j \otimes a_{ji}$, we have
$\varphi_V(v_i^*) = \sum_j \Phi(a_{ij})v_j^*$ for the dual basis $(v_i^*)$
(recall that $\Phi^{-1} = \Phi \circ S$).
It is easy to check that $\varphi_V$ is a map of comodules, because
$\Phi * S = S^{-1} * \Phi$. Also it is easily seen that if
$f : V \longrightarrow W$ is a map of comodules, then
$^tf \circ \varphi_W = \varphi_V \circ f^t$.
In this way we get a bijective natural transformation
$\varphi : {\bf D}_r \longrightarrow {\bf D}_l$.
Finally it is easy to check that $\varphi$ is a morphism of monoidal
functors since $\Phi$ is a character of $A$ and therefore we have
defined a sovereign structure on ${\rm Co}_f(A)$. \square

\begin{coro}
Let $(A,\Phi)$ be a sovereign Hopf algebra. Then there is a sovereign
structure on ${\rm Mod}_f(A)$.
\end{coro}

\noindent
{\bf Proof}. Let $A^0$ be the dual Hopf algebra of $A$. The categories
${\rm Mod}_f(A)$ and ${\rm Co}_f(A^0)$ are monoidally equivalent and
$(A^0, \Phi)$ is a cosovereign Hopf algebra. \square

\medskip

We want a converse of proposition 3.10. We obtain a more precise result in the reconstruction setting: 

\begin{theo}
i) Let \socat be a sovereign monoidal category and let 
$F : \mathcal C \longrightarrow {\rm Vect}_f(k)$ be a monoidal
functor. Then there is a sovereign character on the Hopf algebra
${\rm End}^\vee(F)$.

\noindent
ii) If the monoidal functor $F$ is sovereign (with respect to any
sovereign structure on ${\rm Vect}_f(k)$), then the square of the antipode
of ${\rm End}^\vee(F)$ is equal to the identity.
\end{theo} 

\noindent
{\bf Proof}. i) Let us endow ${\rm Vect}_f(k)$ with the autonomous structure of example 2.9 and with the isomorphism of monoidal functors $\psi : 
{\bf D}_r \longrightarrow {\bf D}_l$ which is the identity.
It is easy to see that $\psi$ is unique for this choice of autonomous
structure. We get in this way a sovereign structure on ${\rm Vect}_f(k)$.
Let us consider the element  $u^{\varphi,\psi} \in {\rm Aut}^\otimes(F)$
of lemma 3.5. By theorem 2.11.iii)  $u^{\varphi,\psi}$ gives rise 
to a character $\Phi$ on ${\rm End}^\vee(F)$.
Let us describe $\Phi$. Let $X \in {\rm ob}(\mathcal C)$ and
$f \in {\rm End}(F(X))$. Then
\begin{eqnarray*}
\Phi([X,f]) & = & {\rm Tr}(f \circ u_X^{\varphi,\psi}) \\
& = & {\rm Tr}(f \circ
\kappa_{F(X)}^{-1} \circ (l_X^{-1}
\circ F(\varphi_X) \circ r_X \circ \psi_{F(X)}^{-1})^t \circ \kappa_{F(X)}) \\
& = & {\rm Tr}(
(l_X^{-1}
\circ F(\varphi_X) \circ r_X \circ \psi_{F(X)}^{-1})^t \circ
\kappa_{F(X)} \circ f \circ \kappa_{F(X)}^{-1}) \\
& = & {\rm Tr}((l_X^{-1}
\circ F(\varphi_X) \circ r_X \circ \psi_{F(X)}^{-1})^t
\circ (^tf)^t) \\
& = & {\rm Tr}(^tf \circ l_X^{-1} \circ F(\varphi_X) \circ r_X
\circ \psi_{F(X)}^{-1})
\end{eqnarray*}
 and  $\Phi^{-1}([X,f]) =  {\rm Tr}
(^tf \circ \psi_{F(X)} \circ r_X^{-1} \circ F(\phi_X)^{-1} \circ l_X).
$
We want to show that $S^{-1}=\Phi * S * \Phi^{-1}$.
Let $\gamma_X = l_X^{-1} \circ F(\varphi_X) \circ r_X
\circ \psi_{F(X)}^{-1}$ and $\beta_X = F(\varphi_X) \circ r_X \circ
\psi_{F(X)}^{-1}$. We have 
\begin{eqnarray*}
S^{-1}([X,f])
 & = & [X^\vee,r_X \circ f^t \circ r_X^{-1}] \\
 & = & [^\vee X, F(\varphi_X) \circ r_X \circ \psi_X^{-1} \circ  \ ^tf \circ
\psi_{F(X)} \circ r_X^{-1} \circ F(\varphi_X^{-1})] \\
& = & [^\vee X,\beta_X \circ \ ^tf \circ \beta_X^{-1}].
\end{eqnarray*}
It is sufficient to work with rank one operators. Let $(e_i)$ be a basis
of $F(X)$  with dual basis $(e_i^*)$ and bidual basis $(e_i^{**})$.
Then we have
$$
S^{-1}([X,e_i^* \otimes e_j]) =
[\ {^\vee X}, \beta_X \circ (e_j^{**} \otimes e_i^*) \circ \beta_X^{-1}]
= [\ {^\vee X}, (e_j^{**} \circ \beta_X^{-1}) \otimes (\beta_X(e_i^*))].
$$ 
On the other hand
\begin{eqnarray*}
& & \Phi * S * \Phi^{-1} ([X, e_i^* \otimes e_j])  \\
& & = \sum_{p,q} {\rm Tr}((e_p^{**} \otimes e_i^*) \circ \gamma_X)
{\rm Tr}((e_j^{**} \otimes e_q^*) \circ \gamma_X^{-1})
[ {^\vee X}, l_X \circ (e_q^{**} \otimes e_p^*) \circ l_X^{-1}] \\
& &  = \sum_{p,q} e_p^{**}(\gamma_X(e_i^*)) e_j^{**}(\gamma_X^{-1}(e_q^*))
[ {^\vee} X, (e_q^{**} \circ l_X^{-1}) \otimes (l_X(e_p^*))] \\
& & = \left[ {^\vee X},
\left(\sum_q e_j^{**}(\beta_X^{-1} \circ l_X(e_q^*)) e_q^{**} \circ l_X^{-1}
\right)\otimes  
\left(\sum_p e_p^{**}(\beta_X \circ l_X^{-1}(e_i^*)) l_X(e_p^*)\right)\right] \\
& & = [{^\vee X}, (e_j^{**} \circ \beta_X^{-1}) \otimes (\beta_X(e_i^*))]
= S^{-1}([X,e_i^* \otimes e_j]).
\end{eqnarray*}
Therefore $S^{-1}=\Phi * S * \Phi^{-1}$ and the proof of i) is complete. 

\noindent
ii) Let us first suppose that $F$ is sovereign with respect
to the sovereign structure defined on ${\rm Vect}_f(k)$ in i).
Then $u^{\varphi,\psi}=1_F$ by lemma 3.5 and $S^{-1} = S$ by i).
The general statement follows now from lemma 3.6 and the fact that the
isomorphism of monoidal functors
 $\psi : {\bf D}_r \longrightarrow {\bf D}_l$ used in the proof of
i) is unique with respect to the choice of the autonomous structure
of example 2.9. \square 
\medskip

\begin{rem}
{\rm There is a unique sovereign structure on ${\rm Vect}_f(k)$
in the sense of \cite{[Ma]}.}
\end{rem}

\begin{rem}
{\rm Theorem 3.12 is still valid if $k$ is a ring (see remark 2.14).}
\end{rem}

\begin{coro}
i) Let $A$ be a Hopf algebra with bijective antipode. Then there is a sovereign character on $A$ if and only of there is a sovereign structure on
${\rm Co}_f(A)$.

\noindent
ii) Let $A$ be a Hopf algebra. Then the square of the antipode of $A$
is equal to the identity if and only if ${\rm Co}_f(A)$ is autonomous
and admits a sovereign structure for which the forgetful functor is
sovereign (with respect to any sovereign structure on ${\rm Vect}_f(k))$.
\end{coro}

\noindent
{\bf Proof}. Statement i) follows from theorem 3.12.i) and theorem 2.12.iii).
Statement ii) follows from theorem 3.12.ii) and theorem 2.12.iii). \square

\medskip

The requirement $S \circ S = id_A$ for a Hopf algebra (such Hopf algebras
are often called involutory) was natural for the purpose of adapting
results from group theory to a more general setting (see for example
 a generalization of Mashke's theorem in \cite{[La]}). However there was
no categorical notion on which this notion could rely.
Corollary 3.15 gives a characterization of involutory Hopf algebras 
from the categorical viewpoint.

\medskip

We end this section with a quick look at dimension theory.
Maltsiniotis has shown in \cite{[Ma]} that a sovereign structure is exactly the one needed to define a trace theory.

\begin{defi}
Let \socat be a sovereign monoidal category and let 
$X \in {\rm ob}(\mathcal C)$. The left (resp. right) dimension of
$X$ is the element of ${\rm End}(X)$ defined by
$${\rm dim}^l_{\varphi}(X) = \varepsilon_X 
\circ (\varphi_X \otimes 1_X) \circ d_X$$
$$({\rm Resp.} \quad  
{\rm dim}^r_{\varphi}(X) = e_X 
\circ (1_X \otimes \varphi_X^{-1}) \circ \eta_X)$$
If $(A,\Phi)$ is a cosovereign Hopf algebra (resp. sovereign Hopf
algebra) the left and right dimension of an object $V$ of
${\rm Co}_f(A)$ (resp. ${\rm Mod}_f(A)$) are denoted by
${\rm dim}^l_{\Phi}(V)$ and ${\rm dim}^r_{\Phi}(V)$.
\end{defi}

Let us note that the left and right dimensions may not coincide.

Let $(A,\Phi)$ be a cosovereign Hopf algebra and let $V$
be a finite-dimensional $A$-comodule with basis $(v_i)$ such that
$\alpha_V(v_i) = \sum_j v_j \otimes a_{ji}$.
Then ${\rm dim}^l_{\Phi}(V) = \sum_i \Phi(a_{ii})$ 
and ${\rm dim}^r_{\Phi}(V) = \sum_i \Phi(S(a_{ii})$.
Let $(A,\Phi)$ be a sovereign Hopf algebra and let $V$
be a finite-dimensional $A$-module. We have
${\rm dim}^l_{\Phi}(V)= {\rm Tr}(\Phi)$
 and ${\rm dim}^r_{\Phi}(V) = {\rm Tr}(\Phi^{-1})$
where $\Phi$ and $\Phi^{-1}$ are considered as operators on $V$. 

\begin{prop}
(see \cite{[Ma]}, corollary 3.5.25).
Let \socat be a sovereign monoi-dal category and let $X$ and $Y$ be
objects of $\mathcal C$.

\noindent
i) ${\rm dim}^l_{\varphi}(I) = 1_I = {\rm dim}^r_{\varphi}(I)$

\noindent
ii) If $X \cong Y$ then ${\rm dim}^l_{\varphi}(X) = {\rm dim}^l_{\varphi}(Y)$ 
and ${\rm dim}^r_{\varphi}(X) = {\rm dim}^r_{\varphi}(Y)$

\noindent
iii) ${\rm dim}^l_{\varphi}( ^\vee X) = {\rm dim}^r_{\varphi}(X) = 
{\rm dim}^l_{\varphi}(X ^\vee)$ ;
${\rm dim}^r_{\varphi}(^\vee X) = {\rm dim}^l_{\varphi}(X) = 
{\rm dim}^r_{\varphi}(X^\vee)$.

\noindent
iv) If ${\rm End}(I)$ is central ($u \otimes f = f \otimes u$ for all
morphisms $f$ of $\mathcal C$ and all $u \in {\rm End}(I)$) then
${\rm dim}^l_{\varphi}(X \otimes Y)) =  
{\rm dim}^l_{\varphi}(X) {\rm dim}^l_{\varphi}(Y)$ 
and ${\rm dim}^r_{\varphi}(X \otimes Y)
= {\rm dim}^r_{\varphi}(X) {\rm dim}^r_{\varphi}(Y)$.
\end{prop}

\section{Universal cosovereign Hopf algebras}

In this section we introduce the universal cosovereign Hopf algebras and 
study some of their properties. By universal we mean that every
finite type cosovereign Hopf algebra is a homomorphic quotient of one of 
them (a quantum subgroup in the language of quantum groups).

\medskip

\noindent
{\bf Notations}. We will use matrix notations. If $u = (u_{ij})$
is a $n \times n$ matrix with values in any algebra,
the transpose matrix of $u$ will be denoted by $^t u$.

\begin{defi}
Let $F \in GL_n(k)$. The algebra $H(F)$ is the universal algebra
with generators $(u_{ij})_{1 \leq i,j \leq n}$,
 $(v_{ij})_{1 \leq i,j \leq n}$ and relations:
$$ u {^tv} = {^tv} u = 1 \quad ; \quad vF {^tu} F^{-1} = F {^tu} F^{-1}v =1$$
\end{defi}

The algebras $H(F)$ are closely related to the universal compact quantum
groups of Van Daele and Wang \cite{[VW]}. We go back to this subject
later in the section.

\begin{theo}
Let $F \in GL_n(k)$. Then $H(F)$ is a Hopf algebra:
\begin{eqnarray*}
{\rm with \ comultiplication} \ \ \Delta(u_{ij}) & = & \sum_k u_{ik} \otimes 
u_{kj} \quad ; \quad \Delta(v_{ij}) = \sum_k v_{ik} \otimes v_{kj} ,\\
{\rm with  \ counit}  \ \ \varepsilon (u_{ij}) & = & \delta_{ij} = 
\varepsilon(v_{ij}), \\
{\rm with \  antipode}  \ \ S(u) & = & {^tv} \quad ; \quad S(v) = F {^tu} F^{-1}. 
\end{eqnarray*}
The antipode of $H(F)$ is bijective and its  inverse is given by
$S^{-1}(u) = {^tF} {^tv} {^tF}^{-1}$ and $S^{-1}(v) = {^tu}$. 
There is a sovereign character $\Phi_F$ on $H(F)$ defined by
$\Phi_F(u) = {^tF}$ and $\Phi_F(v) = F^{-1}$ and hence 
$(H(F), \Phi_F)$ is a cosovereign Hopf algebra.

If $A$ is a Hopf algebra with bijective antipode and if $V$ is a 
finite-dimensional $A$-comodule with coaction
$\alpha_V : V \longrightarrow V \otimes A$ such that
$V^\vee \cong {^{\vee} V}$, then there is a matrix $F \in GL_n(k)$ 
($n = {\rm dim}(V)$), a coaction
$\beta_V : V \longrightarrow V \otimes H(F)$ and a Hopf algebra morphism
$\pi : H(F) \longrightarrow A$ such that
$(1_V \otimes \pi) \circ \beta_V = \alpha_V$.
In particular for every finite type cosovereign Hopf algebra $(A,\Phi)$,
there is a surjective Hopf algebra morphism $\pi : H(F) \longrightarrow A$
for some $F \in GL_n(k)$.
\end{theo}

Theorem 4.2 justify the expression universal cosovereign Hopf algebras
for the algebras $H(F)$.

 \medskip

\noindent
{\bf Proof}. It is easily seen that the maps $\Delta$ and $\varepsilon$
are well defined algebra morphisms and thus $H(F)$ is a bialgebra.
In the same way the formulas of the theorem give rise to a well defined 
anti-homomorphism $S$ which is clearly an antipode for $H(F)$. Once again it is clear that $S$ is bijective with inverse as given in the theorem. 
The character $\Phi_F$ is easily seen to be well defined on $H(F)$.
We have $S^{-1}(u) = {^tF} {^tv} {^tF}^{-1} = \Phi_F(u) S(u) \Phi_F^{-1}(u)$
and $S^{-1}(v) = {^tu} = F^{-1}F {^tu}F^{-1}F = F^{-1} S(v) F =
\Phi_F(v) S(v) \Phi_F^{-1}(v)$.
Therefore $S^{-1} = \Phi_F * S * \Phi_F^{-1}$ and
$(H(F), \Phi_F)$ is a cosovereign Hopf algebra.

Let $A$ be a Hopf algebra with bijective antipode
and let $V$ be a finite dimensional $A$-comodule with basis $v_1,\ldots,v_n$
such that $\alpha_V(v_i) =  \sum_j v_j \otimes a_{ji}$.
We have $V^\vee \cong {^{\vee} V}$ and hence there is a matrix 
$F \in GL_n(k)$ such that $^tS(a)F = F {^tS^{-1}(a)}$
($a$ is the matrix $(a_{ij})$). Then there is a Hopf algebra morphism
$\pi : H(F) \longrightarrow A$ defined by $\pi(u) = a$ and
$\pi(v) = {^tS^(a)}$ (since $(^ta)^{-1} = {^tS^{-1}(a)}$).
A coaction $\beta_V : V \longrightarrow V \otimes H(F)$ is defined 
by $\beta_V(v_i) = \sum_j v_j \otimes u_{ji}$.
It is clear that $\pi$ satisfies the requirement in the theorem.
The last assertion is straightforward. \square

\medskip

The Hopf algebras $H(F)$ clearly show that anything can happen
with the dimension theory of sovereign monoidal categories.
Let $U$ be the obvious $n$-dimensional comodule associated with the 
matrix $u$. Then ${\rm dim}^l_{\Phi_F}(U) = {\rm Tr}(F)$
and ${\rm dim}^r_{\Phi_F}(U) = {\rm Tr}(F^{-1})$. These two scalars
may not coincide.
It is also clear in this example that the dimension of an 
object may take any value with respect to different sovereign structure.

The next result reduces the list of the algebras $H(F)$:

\begin{prop}
Let $F$ and $K \in GL_n(k)$ and let $\lambda \in k^*$.
Then $H(\lambda F) = H(F)$, $H(F) \cong H(KFK^{-1})$ and 
$H(F) \cong H( ^tF^{-1})$ (as Hopf algebras).
\end{prop} 

\noindent
{\bf Proof}. The first statement is clear. A Hopf algebra isomorphism
$\phi : H(F) \rightarrow H(KFK^{-1})$ is defined by
$\phi(u) = {^tK} u {^tK^{-1}}$ and $\phi(v) = K^{-1} v K$.
 A Hopf algebra isomorphism
$\psi : H(F) \longrightarrow H({^tF^{-1}})$ is defined by
$\psi(u) = v$ and $\psi(v) = F u F^{-1}$. \square

\medskip

We now study some cosemisimplicity properties of $H(F)$.
Let us recall that a Hopf algebra is cosemisimple if and only if
there is a Haar measure on it (\cite{[A]}).

\begin{prop}
Let $F \in GL_n(k)$. If ${\rm Tr}(F) = 0$ or ${\rm Tr}(F^{-1}) = 0$
then $H(F)$ is not cosemisimple.
\end{prop}

\noindent
{\bf Proof}. One can assume that the base field is algebraically closed.
The $n$-dimensional comodule $U$ associated to the elements $u_{ij}$
is clearly irreducible. The result follows from
an application to $U$ of Larson's orthogonality relation (\cite{[La]})
as expressed by Woronowicz (\cite{[W1]}, see also \cite{[DK]}, proposition 
3.5). \square

\medskip

We now make contact with the theory of compact quantum groups.
We assume that the reader is familiar with Hopf $*$-algebras 
\cite{[V1]} and with the algebraic theory of compact quantum groups as in
\cite{[DK]}. We now assume the base field to be the field of complex numbers.

\medskip

\noindent
{\bf Notations}. Let $a = (a_{ij})$ where $A$ is a $*$-algebra.
Then the matrix $(a_{ij}^*)$ is denoted by $\overline a$ and 
the matrix $^t \overline a$ is denoted by $a^*$.

\begin{defi}
A Hopf $*$-algebra is a complex algebra $A$, which is a $*$-algebra and 
whose coproduct $\Delta : A \longrightarrow A\otimes A$ is a $*$-homomorphism.
A $CQG$ algebra is Hopf $*$-algebra $A$ such that every finite-dimensional
$A$-comodule is unitarizable.
In other words for every matrix $a=(a_{ij}) \in M_n(A)$ such that
$\Delta(a_{ij}) = \sum_k a_{ik} \otimes a_{kj}$ and
$\varepsilon(a_{ij}) = \delta_{ij}$, there is a matrix $F \in GL_n(\mathbb C)$
such that the matrix $FaF^{-1}$ is unitary.  
\end{defi}

A $CQG$ algebra may be thought as the algebra of representative functions on
a compact quantum group. A Hopf $*$-algebra is $CQG$ if and only if there is
a faithful Haar measure on it (\cite{[DK]}, 3.10). In particular a $CQG$
algebra is cosemisimple. In the next result we find a necessary and 
sufficient condition for the Hopf algebra $H(F)$ to admit a 
$CQG$ algebra structure. We say that a matrix $F \in GL_n(\mathbb C)$ is
relatively positive if there is a scalar $\lambda \in \mathbb C^*$ such that 
$\lambda F$ is a positive matrix.

\begin{prop}
Let $F \in GL_n(\mathbb C)$. Then $H(F)$ admits a $CQG$ algebra structure if
and only if $F$ is conjugate to a relatively positive matrix.
\end{prop}

\noindent
{\bf Proof}. Let us assume that $H(F)$ admits a $CQG$ algebra structure.
Then there is $K \in GL_n(\mathbb C)$ such that $KuK^{-1}$ is unitary
and hence $uK^{-1} K^{-1 *} u^* K^* K = I$ which implies
$^tv = K^{-1} K^{-1 *} u^* K^* K$ and 
$^tu = {^tK} \overline K v^* \overline K^{-1} {^t K}^{-1}$.
There is $Q \in GL_n(\mathbb C)$ such that $QvQ^{-1}$ is unitary and hence
$v Q^{-1} Q^{-1 *} v^* Q^* Q = I$  which implies
$F {^tu} F^{-1} = Q^{-1} Q^{-1 *} v^* Q^* Q$ and 
$^tu = F^{-1} Q^{-1} Q^{-1 *} v^* Q^* Q F$.
We get $Q^* Q F {^tK \overline K} = \alpha I$ for some 
$\alpha \in \mathbb C^*$ (the $v_{ij}$'s are linearly independent).
Therefore $^tK^{-1} F {^tK} = \alpha {^tK}^{-1} Q^{-1} Q^{-1 *} 
\overline K ^{-1}$ and $F$ is conjugate to a relatively positive matrix.

Conversely we can assume that $F$ is a positive matrix by proposition 4.3.
It is easy to see that a Hopf $*$-algebra structure is defined on $H(F)$
by letting $\overline u = v$ (at this point we only use $F^* = F$).
The matrix $u$ is unitary. Let $K = \sqrt F$. It is easy to see that the matrix
$KvK^{-1}$ is unitary.
It follows that $H(F)$ is a $CQG$ algebra since it is generated by the
matrices $u$ and $v$ (\cite{[DK]}, proposition 2.4). \square  

\medskip

When $F$ is a positive matrix, it is easy to see that $H(F)$ is the
$CQG$ algebra of representative functions on the compact quantum group
$A_u(F)$ defined by Van Daele and Wang in theorem 1.3 of \cite{[VW]}
(in fact it is sufficient to consider positive matrices to get the
family of universal compact quantum groups). In that case the representation 
semi-ring of $H(F)$ is described by Banica in \cite{[Ba]}:
the irreducible comodules are labelled by the free product
$\mathbb N * \mathbb N$.

Woronowicz has shown (\cite{[W1]}, theorem 5.6) that a $CQG$ algebra
always has a sovereign character. For a general cosemisimple Hopf 
algebra $A$, there is a convolution invertible $\lambda$
such that $S^2 = \lambda * id * \lambda^{-1}$ (\cite{[La]}) and the
map $\sigma = \lambda * id * \lambda$ is an algebra morphism
(the modular homomorphism of theorem 5.6 in \cite{[W1]}).
However it does not seem to be clear that a cosemisimple Hopf algebra 
always has a sovereign character.
  
\section{Other examples}

We first take a look at finite-dimensional Hopf algebras.

\begin{theo}
Let $(A,\Phi)$ be a finite-dimensional cosovereign (or sovereign)
Hopf algebra over a field of characteristic zero.
If ${\rm dim}_{\Phi}^l(A) \not= 0$ or ${\rm dim}_{\Phi}^r(A) \not= 0$
then $S \circ S = 1_A$ and $A$ is semisimple and cosemisimple.
\end{theo}

\noindent
{\bf Proof}. We can assume that $k$ is algebraically closed and that
$(A,\Phi)$ is a cosovereign Hopf algebra (the sovereign case is dual).
We consider $A$ as an $A$-comodule via the comultiplication.
For every finite-dimensional $A$-comodule there is an $A$-comodule 
isomorphism $A \otimes V \cong A^{{\rm dim}(V)}$ (see proposition 1
in \cite{[B2]}).
By proposition 3.17 we have 
${\rm dim}_{\Phi}^l(A) {\rm dim}_{\Phi}^l(V) = 
{\rm dim}(V) {\rm dim}_{\Phi}^l(A)$ (the dimension is clearly additive
on direct sums). In particular ${\rm dim}_{\Phi}^l(A) = {\rm dim}(A)$
if ${\rm dim}_{\Phi}^l(A) \not= 0$. 
Let $e_1,\ldots,e_n$ be a basis of $A$ such that
$\Delta(e_i) = \sum_j e_j \otimes a_{ji}$.
Let $F$ be the matrix $F = (\Phi(a_{ij}))$. The matrix $F$ can be assumed
to be triangular. We have $\Phi^{*k} = \varepsilon$ for some integer $k$
since $\Phi$ is a character
and hence $F^k = I$. This means that the elements $\Phi(a_{ii})$
are k-th roots of unity.
But ${\rm dim}_{\Phi}^l(A) = {\rm dim}(A) = n = \sum_i \Phi(a_{ii})$
and therefore $\Phi(a_{ii}) = 1$ for all $i$ since 
the base field is of characteristic 0.
This also implies that $F$ is a diagonal matrix and $F=I$.
Then $\Phi = \varepsilon$ and $S\circ S = 1_A$.
Now $A$ is semisimple by \cite{[La]} theorem 4.3
and is also cosemisimple by \cite{[LR]}. The proof is the
same if  ${\rm dim}_{\Phi}^r(A) \not= 0$. \square

\begin{ex}
{\rm It is easy to see that Sweedler's famous 4-dimensional Hopf algebra
(see \cite{[K]}, p. 67)
admits a sovereign element (the group-like element)
and a sovereign character (the only non-trivial character).
In both cases, the left and right dimensions of this algebra are equal to zero.

Sweedler's algebra clearly shows a way to construct sovereign Hopf algebras.
Let $H_n$ be the quotient of the free algebra 
$k\{X_1,\ldots,X_n,\Phi, \Phi^{-1} \}$
by the two-sided ideal generated by the relations
$\Phi \Phi^{-1} = 1  = \Phi^{-1} \Phi$.
Then $H_n$ is a Hopf algebra with comultiplication
$\Delta(X_i) = 1 \otimes X_i + X_i \otimes \Phi$,
$\Delta(\Phi) = \Phi \otimes \Phi$, 
with counit $\varepsilon(X_i) = 0$, $\varepsilon(\Phi) = 1$
and with antipode $S(X_i) = X_i \Phi$ and $S(\Phi) = \Phi^{-1}$.
It is easy to see that $S$ is bijective and that
$\Phi$ is a sovereign element in $H_n$. For more examples
of this kind, see \cite{[Ta],[Ra]}.} 
\end{ex}

We now examine a class of examples closely related to the quantum
groups $SU(n)$ of \cite{[W2]}.

Let $V = k^n$ and let $e_1,\ldots, e_n$ be the canonical basis with
dual basis $e_1^*,\ldots, e_n^*$. Let $N \geq 2$ be an integer and let
$E : V^{\otimes N} \longrightarrow k$ be a linear map.
Let $E(i_1,\ldots,i_N) = E(e_{i_1} \otimes \ldots \otimes e_{i_N})$.
We say that $E$ is left non-degenerate if the linear map
$$V^{\otimes N-1} \longrightarrow V^*, \quad
e_{i_1} \otimes \ldots \otimes e_{i_{N-1}} \longmapsto
\sum_k E(i_1,\ldots,i_{N-1},k) e_k^*$$
is surjective. In this case there are scalars $\lambda(i_1,\ldots,i_N)$
such that
$$ (\star) \quad
\sum_{j_1,\ldots,j_{N-1}} \lambda(i,j_1,\ldots,j_{N-1})
E(j_1,\ldots,j_{N-1},k) = \delta_{ik}, \quad 1\leq i,k \leq n.$$
We say that $E$ is right non-degenerate if the linear map
$$V^{\otimes N-1} \longrightarrow V^*, \quad
e_{i_1} \otimes \ldots \otimes e_{i_{N-1}} \longmapsto
\sum_k E(k,i_1,\ldots,i_{N-1}) e_k^*$$
is surjective. In that case there are scalars $\mu(i_1,\ldots,i_N)$ such that
$$ (\star \star) \quad
\sum_{j_1,\ldots,j_{N-1}} E(k,j_1,\ldots,j_{N-1})
\mu(j_1,\ldots,j_{N-1},i) = \delta_{ik}, \quad 1\leq i,k \leq n.$$

\begin{theo}
Let $E : V^{\otimes N} \longrightarrow k$ be a left and right non-degenerate
linear map. Let $SL(E)$ be the universal algebra with generators
$(a_{ij})_{1\leq i,j\leq n}$ and relations:
$$
(5.3.1) \quad
\sum_{j_1,\ldots,j_N} E(j_1,\ldots,j_N) a_{j_1 i_1} \ldots a_{j_N i_N}
= E(i_1,\ldots, i_N)1, \quad 1\leq i_1, \ldots ,i_N \leq n$$ 
$$
(5.3.2) \quad
\sum_{j_1,\ldots,j_N} E(j_1,\ldots,j_N) a_{i_1 j_1} \ldots a_{i_N j_N}
= E(i_1,\ldots, i_N)1, \quad 1\leq i_1, \ldots, i_N \leq n$$ 
i) Then $SL(E)$ is a Hopf algebra with bijective antipode.

\noindent
ii) Assume that there are invertible scalars $(\beta_i)_{1\leq i \leq n}$ such that
$E(j_1,\ldots,j_{N-1},i) \ = \beta_i E(i,j_1,\ldots, j_{N-1})$ for all
$i,j_1, \ldots, j_{N-1}$. Then there is a sovereign character
$\Phi_{\beta}$ on $SL(E)$ such that $\Phi_{\beta}(a_{ij}) = \delta_{ij}
\beta_i$.

\noindent
iii) If $k$ is a field of characteristic zero and if the field
$\mathbb Q(E(j_1,\ldots, j_N)_{1\leq i_1, \ldots, i_N \leq n})$
can be ordered, then $SL(E)$ is cosemisimple.

\noindent
iv) If $k = \mathbb C$ and $E(j_1, \ldots, j_N) \in \mathbb R$
for all $j_1, \ldots, j_N$, then $SL(E)$ admits a $CQG$ algebra structure. 
\end{theo}

\noindent
{\bf Proof}. i) It is easily seen that $SL(E)$ is a bialgebra with coproduct
$\Delta(a_{ij}) = \sum_k a_{ik} \otimes a_{kj}$ and counit $\varepsilon(a_{ij})
= \delta_{ij}$. Let us show that the matrix $a = (a_{ij})$ is invertible.
Let us consider equation 5.3.1. Multiplying  by
$\lambda(k,i_1,\ldots,i_{N-1})$ and summing over $i_1, \ldots, i_{N-1}$,
we get that $a$ is left invertible (we use ($\star$)). In the same way
$a$ is right invertible (use 5.3.2 and ($\star \star)$) and 
therefore $a$ is invertible and $SL(E)$ is a Hopf algebra
by \cite{[Wa]}, theorem 1. Let us show that $^ta$ is invertible.
By 5.3.2 and ($\star$) $^ta$ is right invertible and by
5.3.1 and ($\star \star$) $^ta$ is left invertible.
Hence the antipode of $SL(E)$ is invertible.

\noindent
ii) The character $\Phi_{\beta}$ is easily seen to be well defined.
Let us consider equation 5.3.2: multiplying on the left by
$S(a_{k i_1})$ and summing over $i_1$, we get
$$
(\star \star \star) \quad
\sum_{j_2,\ldots,j_N} E(k,j_2,\ldots j_N) a_{i_2 j_2} \ldots a_{i_N j_N}
= \sum_{i_1} S(a_{k i_1}) E(i_1,\ldots, i_N). $$ 
We have
\begin{eqnarray*}
& & \sum_i \sum_k E(k,i_2,\ldots,i_N) S(a_{ik}) \beta_i a_{ji} \\
& = & \sum_i \sum_{j2,\ldots,j_N} E(i,j_2,\ldots,j_N) \beta_i
a_{i_2j_2} \ldots a_{i_N j_N} a_{ji} \quad {\rm by} \ (\star \star \star) \\
& = & \sum_i \sum_{j2,\ldots,j_N}
E(j_2,\ldots j_N,i) a_{i_2j_2} \ldots a_{i_N j_N} a_{ji} \\
& = & E(i_2,\ldots,i_N,j)\quad ({\rm by \ 5.3.2}) = \beta_j E(j,i_2,\ldots,i_N).
\end{eqnarray*}
Using ($\star \star$) we get
$$\sum_i S(a_{il}) \beta_i a_{ji} = \delta_{jl} \beta_j 
\ {\rm for \ all} \ j,l.$$
The inverse of the matrix $^ta$ is $^tS^{-1}(a)$ and hence we have
$S^{-1}(a_{il}) = \beta_i S(a_{il}) \beta_j^{-1} = 
\Phi_{\beta} * S * \Phi_{\beta}^{-1}(a_{il})$.
This means that $\Phi_{\beta}$ is a sovereign character on $SL(E)$.

\noindent
iii) There is an algebra automorphism $\tau$ of $SL(E)$ defined by
$\tau(a_{ij}) = a_{ji}$. Hence by \cite{[Bi]}, 4.7 $SL(E)$ is cosemisimple.
Statement iv) follows from \cite{[Bi]}, 4.6. \square

\medskip 

In a special case the $SL(E)$ construction leads to the quantum groups $SL_q$.
Let $N=n$ and let $q \in k^*$. Let 
$E_q : V^{\otimes n} \longrightarrow k$ defined by
$E_q(i_1, \ldots i_N) = 0$ if two indices are equal and otherwise
$E_q(i_1, \ldots i_N) = (-q)^{l(\sigma)}$ where $l(\sigma)$ is the length
of the permutation $\sigma(k) = i_k$.
The Hopf algebras $SL(E_q)$ and $SL_q(n)$ are isomorphic.
This fact can be proved using the same proof as in Rosso's comparison
of the quantum $SL$ groups of Woronowicz (\cite{[W2]})
and Drinfeld (\cite{[Dr]})
(\cite{[Ro]}, theorem 6). See also \cite{[V1]} for useful computations.
The elements $\beta_i$ of theorem 5.3 are given by
$\beta_i = (-q)^{n+1-2i}$ and therefore  $SL_q(n)$ admits a sovereign character
$\Phi$ defined by $\Phi(a_{ij}) = \delta_{ij} (-q)^{n+1-2i}$.

\appendix
\section{Appendix}

In this appendix we explicitly write and prove the correspondence between
cotwists and sovereign characters for a cobraided Hopf algebra.
We intensively use Sweedler's notations \cite{[Sw]} :
if $a$ is an element of a Hopf algebra we write
$\Delta(a) = \sum a_1 \otimes a_2$.
We first recall some basic definitions (see \cite{[Doi],[K],[JS1]}).
 
\begin{defi}
A cobraided Hopf algebra is a pair $(A, \sigma)$ where $A$ is a Hopf algebra
and $\sigma : A \otimes A \longrightarrow k$ is a convolution invertible
linear map satisfying:

\noindent
(A.1) $\sigma * m = m^{\rm op} * \sigma$, ie for all $x, y \in A$ we
have : 
$\sum \sigma(x_1,y_1) x_2 y_2 =\sum y_1 x_1 \sigma(x_2,y_2)$.

\noindent
(A.2) $\sigma(xy,z) = \sum \sigma(x,z_1) \sigma(y,z_2)$
 for all $x, y, z \in A$.

\noindent
(A.3) $\sigma(x,yz) = \sum \sigma(x_1,z) \sigma(x_2,y)$
 for all $x, y, z \in A$.

A cotwist on a cobraided Hopf algebra $(A, \sigma)$ is a central convolution
invertible linear form $\tau$ on $A$ satisfying
$\tau \circ m = {^t \sigma} * (\tau \otimes \tau) * \sigma$
(where ${^t \sigma}(x,y) = \sigma(y,x)$), ie for all $x,y \in A$
we have $\tau(xy) = \sum \sigma(y_1, x_1) \tau(x_2) \tau(y_2)
\sigma(x_3,y_3)$.
\end{defi}

Let $\sigma^{-1}$ be the convolution inverse of $\sigma$. The following 
equalities hold for all $x, y, z \in A$ (see \cite{[Doi]}):

\noindent
(A'.1) $\sum \sigma^{-1}(x_1,y_1) y_2 x_2 =\sum x_1 y_1 \sigma^{-1}(x_2,y_2)$.

\noindent
(A'.2) $\sigma^{-1}(xy,z) = \sum \sigma^{-1}(y,z_1) \sigma^{-1}(x,z_2)$

\noindent
(A'.3) $\sigma^{-1}(x,yz) = \sum \sigma^{-1}(x_1,y) \sigma^{-1}(x_2,z)$

\noindent
(A.4) $\sigma(1,x) = \varepsilon(x) = \sigma(x,1)$ ; 
 \ (A'.4) $\sigma^{-1}(1,x) = \varepsilon(x) = \sigma^{-1}(x,1)$

\noindent
(A.5) $\sigma^{-1}(x,y) = \sigma(S(x),y)$ ; 
(A.6) $\sigma(x,y) = \sigma^{-1}(x,S(y))$ ; 
(A.7) $\sigma(x,y) = \sigma(S(x),S(y))$.

\medskip

Let $\lambda$ be the linear form on $A$ defined by
$\lambda(x) = \sum \sigma(x_1,S(x_2))$. Then $\lambda$ is invertible with 
inverse given by $\beta(x) = \sum \sigma^{-1}(S(x_1),x_2)$. Furthermore
$S^2 = \beta * id * \lambda$ (\cite{[Doi]}, theorem 1.3) and in particular 
the antipode of a cobraided Hopf algebra is bijective.

We first observe the following result: 

\begin{lemm}
Let $(A,\sigma)$ be a cobraided Hopf algebra. Then $\beta$ satisfies the 
cotwist equation $\beta \circ m = {^t \sigma} * (\beta \otimes \beta)
* \sigma$.
\end{lemm}

\noindent
{\bf Proof}. Let $x,y \in A$. We have
\begin{eqnarray*}
& & {^t \sigma}^{-1} * (\beta \circ m) * \sigma^{-1}(xy) \\
& = & 
\sum \sigma^{-1}(y_1,x_1) \sigma^{-1}(S(x_2y_2),x_3y_3)
\sigma^{-1}(x_4,y_4) \\
& = &
\sum \sigma^{-1}(y_2,x_2) \sigma^{-1}(S(y_1x_1),x_3y_3)
\sigma^{-1}(x_4,y_4) \quad ({\rm by} \ \sigma^{-1}*m^{\rm op}
= m * \sigma^{-1}) \\
& = &
\sum \sigma^{-1}(y_2,x_2) \sigma^{-1}(S(y_1),x_3y_3)
\sigma^{-1}(S(x_1),x_4y_4) \sigma^{-1}(x_5,y_5) \quad ({\rm by \ (A'.2)}) \\
& = &
\sum \sigma^{-1}(y_3,x_2) \sigma^{-1}(S(y_2),x_3)
\sigma^{-1}(S(y_1),y_4) \sigma^{-1}(S(x_1),x_4y_5) \sigma^{-1}(x_5,y_6)
 ({\rm by \ (A'.3)}) \\
& = &
\sum \sigma^{-1}(S(y_1),y_2) \sigma^{-1}(S(x_1),x_2 y_3)
\sigma^{-1}(x_3,y_4) \quad ({\rm by \ (A'.2) \ and \ (A'.4)}) \\
& = &
\sum \beta(y_1) \sigma^{-1}(S(x_1), x_2y_2) \sigma^{-1}(x_3,y_3) \\
& = &
\sum \beta(y_1) \sigma^{-1}(S(x_1),y_3x_3) \sigma^{-1}(x_2,y_2)
\quad ({\rm by} \ m * \sigma^{-1} = \sigma^{-1} * m^{\rm op}) \\
& = &
\sum \beta(y_1) \sigma^{-1}(S(x_2),y_3)
\sigma^{-1}(S(x_1),x_4) \sigma^{-1}(x_3,y_2) \quad ({\rm by \ (A'.3)}) \\
& = &
\sum \beta(y_1) \sigma^{-1}(S(x_2)x_3,y_2) \sigma^{-1}(S(x_1),x_4)
\quad ({\rm by \ (A'.2)}) \\
& = &
\beta(x) \beta(y) \quad {\rm by \ (A'.4).} \ \square    
\end{eqnarray*}

On the Hopf algebra level, proposition 2.11 from \cite{[Ye]} takes
the following form: 

\begin{theo}
Let $(A,\sigma)$ be a cobraided Hopf algebra. There is a bijective 
correspondence between sovereign characters on $A$ and cotwists on
$(A,\sigma)$. Explicitly we have: 

\noindent
i) If $\Phi$ is a sovereign character on $A$ then 
$\tau = \Phi * \beta$ is a cotwist on $(A, \sigma)$.

\noindent
ii) If $\tau$ is a cotwist on $(A,\sigma)$ then $\Phi = \tau * \beta^{-1}$
is a sovereign character on $A$ 
\end{theo}

\noindent
{\bf Proof}. i) Let us first show that $\tau$ is central. For this
purpose it is sufficient to prove that $id* \tau = \tau * id$.
We have $id * \tau = id * \Phi * \beta = \Phi * S^2 * \beta = 
\Phi * \beta * id = \tau * id$ since $S^2 * \beta = \beta * id$ and 
$\Phi * S^2 = id * \Phi$. Let $x,y \in A$. We have
\begin{eqnarray*}
\tau(xy) & = & \sum \Phi(x_1) \Phi(y_1) \beta(x_2y_2)
\quad ({\rm since} \ {\rm \Phi \ is \ a \ character}) \\
& = & 
\sum \Phi(x_1) \Phi(y_1) \sigma(y_2,x_2)
\beta(x_3) \beta(y_3) \sigma(x_4,y_4) \quad ({\rm by \ lemma \ A.2}) \\
& = &
\sum \sigma(S^{-2}(y_1),S^{-2}(x_1)) \Phi(x_2) \Phi(y_2)
\beta(x_3) \beta(y_3) \sigma(x_4,y_4) \quad 
(\Phi * id = S^{-2} * \Phi) \\
& = &
\sum \sigma(y_1,x_1) \Phi(x_2) \beta(x_3) \Phi(y_2) \beta(y_3)
\sigma(x_4,y_4) \quad ({\rm by \ (A.7)}) \\
& = &
{^t \sigma} * (\tau \otimes \tau) * \sigma(x,y).
\end{eqnarray*}
Therefore $\tau$ is a cotwist on $(A,\sigma)$.  

\noindent
ii) Let us first remark that 
$\beta^{-1}(x) = \sum \sigma(x_1,S(x_2)) = \sum \sigma(S^{-1}(x_1),x_2)$
by (A.7) and that $\Phi = \beta^{-1} * \tau$ since $\tau$ is central.
Let $x,y \in A$. We have 
\begin{eqnarray*}
\Phi(xy) & = & \sum \beta^{-1}(x_1y_1) \tau(x_2y_2) \\
& = & 
\sum \sigma(S^{-1}(y_1) S^{-1}(x_1),x_2y_2) \tau(x_3y_3) \\
& = &
\sum \sigma(S^{-1}(y_1),x_2y_2) \sigma(S^{-1}(x_1),x_3y_3)
\tau(x_4y_4) \quad ({\rm by \ (A.2))} \\
& = &
\sum \sigma(S^{-1}(y_2),y_3) \sigma(S^{-1}(y_1),x_2)
\sigma(S^{-1}(x_1),x_3y_4) \tau(x_4y_5) \quad ({\rm by \ (A.3)})\\
& = &
\sum \beta^{-1}(y_2) \sigma(S^{-1}(y_1),x_2) \tau(x_3y_3)
\sigma(S^{-1}(x_1),x_4y_4) \quad (\tau \ {\rm central}) \\
& = &
\sum \beta^{-1}(y_1) \sigma(S(y_2),x_2) \tau(x_3y_3)
\sigma(S^{-1}(x_1),x_4y_4) \quad (S^{-1}*\beta^{-1} = \beta^{-1}*S) \\
& = &
\sum \beta^{-1}(y_1) \sigma^{-1}(y_2,x_2) \tau(x_3y_3)
\sigma(S^{-1}(x_1),x_4y_4) \quad ({\rm by \ (A.5))} \\
& = &
\sum \beta^{-1}(y_1) \tau(x_2) \tau(y_2) \sigma(x_3,y_3)
\sigma(S^{-1}(x_1),x_4y_4) \quad (\tau \ {\rm twist}) \\
& = &
\sum \beta^{-1}(y_1) \tau(x_3) \tau(y_2) \sigma(x_4,y_3)
\sigma(S^{-1}(x_2),y_4) \sigma(S^{-1}(x_1),x_5) \quad ({\rm by \ (A.3)})\\
& = &
\sum \beta^{-1}(y_1) \tau(x_4) \tau(y_2) \sigma(x_3,y_3)
\sigma(S^{-1}(x_2),y_4) \sigma(S^{-1}(x_1),x_5) 
\quad (\tau \ {\rm central}) \\
& = &
\sum \beta^{-1}(y_1) \tau(x_4) \tau(y_2) \sigma(x_3 S^{-1}(x_2),y_3)
\sigma(S^{-1}(x_1),x_5) \quad ({\rm by \ (A.2)}) \\
& = &
\sum \beta^{-1}(y_1) \tau(x_2) \tau(y_2) \sigma(S^{-1}(x_1),x_3)
\quad ({\rm by \ (A.4)}) \\
& = &
\sum \beta^{-1}(y_1) \tau(y_2) \tau(x_3) \sigma(S^{-1}(x_1),x_2)
\quad (\tau \ {\rm central}) \\
& = &
\sum \beta^{-1}(y_1) \tau(y_2) \beta^{-1}(x_1) \tau(x_2) = 
\Phi(x) \Phi(y).
\end{eqnarray*}
Therefore $\Phi$ is a character.
Now $\Phi * S * \Phi^{-1} = \beta^{-1} * \tau * S * \tau^{-1} * \beta
= \beta^{-1} * S * \beta = S^{-1}$ ($\tau$ is central):
$\Phi$ is a sovereign character. \square

\bigskip

D\'epartement des Sciences Math\'ematiques,
case 051

Universit\'e Montpellier II

Place Eug\`ene Bataillon, 34095 Montpellier Cedex 5

{\tt e-mail : bichon\char64math.univ-montp2.fr}

\end{document}